\documentclass[a4paper,8pt]{article}
\usepackage{amsmath}
\usepackage{amsthm,amsfonts,amssymb}
\usepackage[dvips]{graphicx}
\usepackage[english]{babel}
\usepackage{color}
\usepackage{graphics}
\usepackage{graphpap}

\providecommand{\U}[1]{\protect\rule{.1in}{.1in}}

\providecommand{\U}[1]{\protect\rule{.1in}{.1in}}
\newtheorem{theorem}{Theorem}

\newtheorem{lemma}[theorem]{Lemma}
\newtheorem{proposition}[theorem]{Proposition}

\makeatletter
\def\blfootnote{\xdef\@thefnmark{}\@footnotetext}
\makeatother
\newcommand{\à}{\`{a}}

\begin{document}

\title{A stochastic model for distributed real time streaming Evolution of the network structure}
\author{Andrea Monsellato}
\maketitle

\begin{abstract}
This paper is a short summary of the main results in the thesis \cite{Mon}.\\
Based on the P2P paradigm we construct a stochastic model for a live media streaming content delivery  network. Starting from the behavior of the out degree process of each node we analyze the evolution of the resulting random graph.
For the out degree process we provide an iterative formula for the distribution probabilities and calculate an explicit expression for the transient distribution.
Moreover we discuss the steady state distribution and the problem of local disconnections.
\end{abstract}

\noindent
\section{Introduction}

Multimedia applications have experienced explosive growth in the last decade and have become a pervasive application over the Internet. A video streaming system distributes a video stream synchronously to a large population of viewers. Today, one-to-many commercial streaming solutions are offered via content delivery networks (CDNs) such as Akamai, Limelight or VitalStream. Such content delivery networks are based on P2P (\emph{peer to peer}) paradigm, i.e. all users (peers) contribute on distribution service.\\
In Section $2$  we  present  a model for a live media streaming content delivery network. The construction of the model is based on P2P paradigm according on \cite{Fav},\cite{Fav1}.
We consider that the number of nodes which lie in the network is managed by a stochastic process and give a stochastic rule to create connections among the nodes. We focus our analysis on the study of the evolution of the structure starting from the local behavior of in and out-degree of a single node (user).\\
In Section $3$ we discuss the local characteristics of the network in the case of zero mortality.
In this case the process associated to the out-degree is a  pure birth process. After having provided the distribution probabilities of this process and its subordinated chain, we show a relationship with a Bernoullian scheme for which the successes probabilities are dependent by the number of previous successes. This approach is very useful for recovering some interesting properties like mean, variance and for recovering concentration and deviation bounds.\\
In Section $4$ we turn to the general case. In this case the evolution of the out-degree process of the generic node is equivalent to a problem in queuing theory, this process is known in literature as discouragement queue. Although from the theoretical point of view the problem of the determination of the transition functions has been solved \cite{Karlindbdp}, the explicit form of them is not present in literature. Assuming that the solution is represented by a Taylor series we obtain an exact formulation of the distribution probabilities being determined recursively, then we demonstrate the convergence of the series. Moreover for the related embedded chain we obtain an explicit expression of the transition functions. Furthermore we discuss the out-degree distribution of the generic node in the steady state and analyze a mechanism for recovering the connections in the case of local disconnections. We propose a solution in terms of "extra job" for the root (master server) and give an estimate for the mean of this extra job.

\section{The Model}

We consider a Poisson process, with parameter $\lambda>0$, that counts the number of accesses to the network. Any access corresponds
to one user who wants to benefit of the streaming service. We consider that the life-time of each user is exponentially distributed with parameter $\mu>0$.
Thus the number of nodes that lie in the network at some time $t>0$ is subordinated to an $M/M/\infty$ queue.\\
Remembering that in the P2P context each node in the network participates in distributing the service to a subset of (\emph{future}) nodes, we give the following rule for generating connections among the nodes:\\

\emph{At time of the access each node sends a service request concurrently to all nodes already in the system.
Any request will be accepted with probability $\frac{1}{1+k}$, where $k$ is the out-degree of the node that receives the request.}\\

\noindent
For each accepted request we generate a directed connection (\emph{directed arc}) from server node to client node.\\
The connection probabilities have been chosen considering that every node decreases its capacity (\emph{availability to service}) while the number of its clients increases.\\
When a node drops from the network we erase this node together all its arcs.\\
We consider the first node as the root (\emph{master server}) and suppose this node to have an infinite life-time.
We observe that by construction the out-degree processes of the nodes are mutually independent and
the number of \emph{in and out} connections of a node is a stochastic process.\\
The generated graph is directed and acyclic, in fact all arcs are directed so that no one of the nodes receives informations
from the nodes which are arrived later.

\subsection{Recovering out-degree processes}

Let $\{N(t)\}_{t\geq 0}$ be the process that counts the number of nodes, except the root, in the networks.
We know that $N(t)$ is, in particular, a birth and death process with parameters $\lambda_l=\lambda$, $l\geq 0$, and $\mu_l=\mu l$, $l\geq 1$ and $\mu_0=0$.\\
Let us consider a generic node that is in the network at time $t>0$. From the local point of view also the out-degree process of the generic node is managed
by a birth and death process.\\
Since the birth rates depend only on the parameter of Poisson process and on the probability of acceptance of the request, we have $\lambda_k=\frac{\lambda}{1+k}$, $k\geq 0$.\\
As regards the mortality rates, let us observe that if the process $N(t)$ has a backward transition the state of the system changes because one of the nodes will be erased together with all its arcs.
This event corresponds to a backward transition of the out-degree of the nodes that were adjacent to the node which has been erased.\\
Assuming that there are $n$ nodes in the network at time $t$, for a generic out-degree process the mortality rate will be the product of the mortality rate of the process $N(t)$, i.e. $n\mu$, and the proportion of the number of neighbors of the node in relation to the total number of nodes.Thus if the generic out-degree was in state $k\leq n$ then $\mu_k=\frac{k}{n}n\mu=k\mu$.\\
In conclusion being $\{X(t)\}_{t\geq 0}$ the generic out-degree process, let $P(X(t)=k)\doteq p(k,t)$ then for all $k\in N_0$ we have:

\begin{equation}\label{equ:outdegreebdp}
p'(k,t)=-\left(\frac{\lambda}{1+k}+k\mu\right) p(k,t)+ \frac{\lambda}{k}p(k-1,t)+1_{\{k\geq 1\}}\mu(k+1)p(k+1,t)
\end{equation}

\noindent
with initial conditions $p(k,0)=\delta_{k,0}$.\\

\noindent
In the following we discuss separately the case $\mu=0$ and $\mu>0$.

\section{Case $\mu=0$}

In this case the equation (\ref{equ:outdegreebdp}) is reduced to

\begin{equation}\label{equ:bc}
\left\{
\begin{array}{ll}
p'(0,t)=\lambda p(0,t)  \\
p'(k,t)= \frac{\lambda}{k}p(k-1,t)-\frac{\lambda}{1+k}p(k,t),\quad if\quad k\neq 0
\end{array}
\right.
\end{equation}

\noindent
with initial conditions $p(k,0)=\delta_{k,0}$.\\

\noindent
A direct use of the Laplace transform gives us the following result.

\begin{proposition}

Let

\begin{align}\label{solbcont}
p(k,t)=\frac{1}{k!}\sum_{j=1}^{k+1}(-1)^{k+1-j}j^k \binom{k+1}{j} e^{-\frac{\lambda t}{j}}
\end{align}

then Eq.(\ref{solbcont}) is the solution of Eq.(\ref{equ:bc}) with initial condition $p(k,0)=\delta_{k,0}$.

\end{proposition}

\noindent
Many interesting properties of the network do not depend on the elapsed time between nodes arrival, i.e. the Poisson process which counts the number of nodes in the network at generic time $t$, but only on the number of the nodes. From this point of view if we consider the subordinated chain $X:N_0 \rightarrow N_0$ of the generic out-degree process with transition probabilities $p_{k,k+1}=\frac{1}{1+k}=1-p_{k,k}$, $k\geq 0$, using generating function approach we have the following

\begin{theorem}

Let $\{X_n\}_{n\geq 0}$ the subordinated chain of the birth process $\{X(t)\}_{t\geq 0}$ associated to Eq.(\ref{equ:bc}) with $\lambda=1$, if we put $P(X_n=k)=p_{n,k}$ then

\begin{align*}
\left\{
\begin{array}{ll}
p_{n+1,k}=\frac{k}{k+1}p_{n,k}+\frac{1}{k}p_{n,k-1}, \quad\forall n\geq 0, \quad 1\leq k \leq n+1\\
p_{i,j}=0,\quad i\geq 0,\quad j\geq i+1\\
p_{0,0}=1,\quad p_{n,0}=0,\quad n\geq 1
\end{array}
\right.
\end{align*}

Furthermore

\begin{equation}\label{equ:pnkembchain}
p_{n,k}=\frac{1}{k!}\sum_{i=1}^{k}(-1)^{k-i}\binom{k+1}{i+1}i^{k}\left(\frac{i}{i+1}\right)^{n-k},\phantom{x}n\geq k\geq1
\end{equation}

\end{theorem}

\noindent
As a consequence of Eq.(\ref{equ:pnkembchain}) we establish a useful Bernoullian scheme for the subordinated chain.

\begin{proposition}

Let $Z_n=\sum_{i=0}^{n}Y_i$, $n\geq 0$, where

\begin{equation*}
Y_n|Y_0,...,Y_{n-1}\sim Ber\left(\frac{1}{1+\sum_{j=0}^{n-1}Y_j}\right);\quad Y_0=0
\end{equation*}

then $Z_n\stackrel{d}{=}X_n$, where the probability distribution of $X_n$ is given by Eq.(\ref{equ:pnkembchain}).\\

\end{proposition}

\noindent
The above Bernoullian scheme can be interpreted as a special type of Polya urn:\\

\emph{Consider an urn containing only \textbf{one} white ball and an arbitrary number of red balls.
If we draw the white ball then we add a red ball, while if we draw a red ball we do not do anything.
In both cases we reinsert the drawn ball and proceed to the next drawing.}\\

\noindent
From the previous proposition and some results due to MacDiarmid \cite{McDiarmidcentering} we establish the following deviation bounds.

\begin{proposition}

Let $\{X_n\}_{n\geq 0}$ the process with probability distribution function as in Eq.(\ref{equ:pnkembchain}), it results that

\begin{equation*}
E(X_n)=\frac{-1+\sqrt{1+8n}}{2}
\end{equation*}

Furthermore

\begin{align*}
&P(X_n\geq(1+\epsilon)E(X_n))\leq \exp\left(-\frac{1}{3}\epsilon^2 E(X_n)\right) \quad 0<\epsilon <1;\\&
P(X_n\leq(1-\epsilon)E(X_n))\leq \exp\left(-\frac{1}{3}\epsilon^2 E(X_n)\right) \quad 0< \epsilon <1
\end{align*}

\end{proposition}

\subsection{The in-degree process: mean and deviation bound}

\noindent
Let us remember that when a node have accessed to the network it sends a request to the nodes already in the system to establish a connection.
The accepted requests are subordinated to out-degree processes because any possible transition for a generic out-degree process depends on its associated urn drawing.

\noindent
Let us consider $n+1-i$ trials, $i\leq n$, of the $i$-th urn associated with the corresponding out-degree process.\\
\noindent
If $X^-_{n+1}$ is the in-degree process of the generic $(n+1)$-th node then

\begin{equation}\label{equ:dindegree}
X^-_{n+1}=\sum_{i=1}^{n}Y^i_{n+1}
\end{equation}

\noindent
We observe that this value will not change after the $(n+1)$-th time step, and remember that the event $\{Y^i_{n+1}=1\}$ corresponds to have generated the arc $i\rightarrow n+1$. Moreover being the urn schemes mutually independent then $X^-_{n+1}$ in Eq.(\ref{equ:dindegree}) is a sum of independent Bernoullian r.v.'s. Now we can
give an estimation of the mean of the generic in-degree process.

\begin{lemma}

Let $X^-_{n+1}$ as in Eq.(\ref{equ:dindegree}) then

\begin{equation}\label{equ:indegreemean}
E(X^-_{n+1})\sim \frac{2\sqrt{2}}{3}\sqrt{n}
\end{equation}

\end{lemma}

\noindent
Before giving a deviation bound for the in-degree process we need to estimate the probabilities of success for the $n+1-i$, $n\geq i$,
drawing of the $i$-th urn, i.e. of the event $\{Y^i_{n+1}=1\}$.

\begin{lemma}
\noindent
Let $\theta^i_{n+1}=P(Y^i_{n+1}=1)$, with $i=1,...,n$ then

\begin{align*}
&\theta^i_{n+1}\sim \frac{1}{\sqrt{2(n-i)+1}}
\end{align*}

\end{lemma}

\noindent
Using an approach based on \emph{self-bounding} functions, see \cite{McDiarmidtalagrand}, we establish a two sides deviation bound for the generic in-degree process.

\begin{proposition}

Let $X^-_{n+1}$ as in Eq.(\ref{equ:dindegree}) then the sequence $X^-_{n+1}$ is $(1,0)$-self-bounding and we have

\begin{align*}
P(X^-_i-E(X^-_i)\geq \epsilon E(X^-_i))\leq \exp\left(-\frac{\epsilon^2 E(X^-_i)}{4}\right)\sim\exp\left(-\frac{\epsilon^2 \sqrt{2i}}{4}\right)
\end{align*}

\begin{align*}
P(X^-_i-E(X^-_i)\leq -\epsilon E(X^-_i))\leq \exp\left(-\frac{3\epsilon^2 E(X^-_i)}{8}\right)
\sim\exp\left(-\frac{3\epsilon^2 \sqrt{2i}}{8}\right)
\end{align*}

\end{proposition}

\subsection{Spanning trees}

\noindent
Starting from the in-degree process we establish a function that counts the number of spanning trees and recovers the mean.\\
We know that the arcs $i\rightarrow i+1$, with $i\in N$, exist with probability $1$ (by construction) and there is at least
a spanning tree for any generated graph.\\
\noindent
Let $G_{n+1}$ be the graph at time step $n+1$ and let $T_{n+1}$ be the number of spanning trees in $G_{n+1}$.
Recalling that at time step $n+1$ we add the node $n+1$ and, with their probabilities, the arcs $i\rightarrow n+1$, with $i=1,...,n$,
the number of these generated arcs represents the in-degree of the $n+1$-th node.\\
Observing that

\begin{equation*}
T_{n+1}=\sum_{i=1}^{n}T^i_{n+1}
\end{equation*}

\noindent
where $T^i_{n+1}=1_{(i,n+1)}T_n$, being $1_{(i,n+1)}$ the indicator function of the arc $i\rightarrow n+1$, i.e. $T^i_{n+1}$ is the number of spanning trees which contain the arc $i\rightarrow n+1$, if it has been generated, then from Eq.(\ref{equ:dindegree}) we have

\begin{align*}
T_{n+1}=\prod_{j=1}^{n}\sum_{i=1}^{j}Y^i_{j+1}
\end{align*}

\noindent
Because for every fixed $i$ the $Y^i_{n}$, with $n=i,i+1,...$ are negatively correlated Bernoulli r.v.'s and for every
fixed $n\geq 1$ the $Y^i_{n}$, with $i\leq n-1$, are mutually independent Bernoulli r.v.'s, from Eq.(\ref{equ:indegreemean})
computing the mean of the number of spanning trees in $G_{n+1}$ we have

\begin{align*}
&E(T_{n+1})=\prod_{j=1}^{n}E\left(\sum_{i=1}^{j}Y^i_{j+1}\right)\sim\left(\frac{8}{9}\right)^{\frac{n}{2}}\sqrt{n!}
\end{align*}

\section{Case $\mu>0$}

\noindent
If $\mu>0$ the out-degree process of the generic node is a so called discouragement queue, see \cite{Conollysdq}.
In the following theorem we provide an exact formulation of its distribution probabilities being determined recursively.

\begin{theorem}
Let

\begin{align*}
\left\{
\begin{array}{ll}
S^{(1)}_i=S^{(0)}_{i+1}+b_0 S^{(0)}_i\\
S^{(k+1)}_i=S^{(k)}_{i+1}+b_k S^{(k)}_i-\alpha^2 S^{(k-1)}_i,\quad \forall i\geq 0, \quad 0\leq k\leq i+1\\
b_k=\frac{1}{k+1}+k\alpha^2 \quad k\geq 0, \quad \alpha\doteq\sqrt{\frac{\mu}{\lambda}}
\end{array}
\right.
\end{align*}

being $S^{(-1)}_i=0 \quad \forall i\geq 0$ and $S^{(n)}_i=0 \quad \forall n>i$, then $\forall k\geq 0$

\begin{align}\label{seriesbd}
p(k,t)=\frac{1}{\alpha^{2k} k!}\sum_{i=0}^{+\infty}\frac{(\lambda t)^i}{i!}S^{(k)}_{i},\quad t\geq 0
\end{align}

is the solution of Eq.(\ref{equ:outdegreebdp}), with initial conditions $p(k,0)=\delta_{k,0}$.

\end{theorem}

\noindent
The above result is obtained making use of Taylor series and showing that the coefficient $S^{(k)}_{i}$ are bounded by Bessel numbers. For more details on Bessel numbers see \cite{bessel}.\\

\noindent
As regards the probability distribution of the embedded chain of the process $\{X(t)\}_{t \geq 0}$ associated to Eq.(\ref{equ:outdegreebdp}), proceeding by induction we establish the following result.

\begin{proposition}

Let $(p_{n,k})_{n,k\geq 0}$ the embedded process of the birth and death process $\{X(t)\}_{t \geq 0}$ associated to Eq.(\ref{equ:outdegreebdp})
it result that

\begin{align*}
\left\{
\begin{array}{ll}
p_{n+1,k}=\frac{1}{1+k(k-1)\alpha^2}p_{n,k-1}+\frac{(k+2)(k+1)\alpha^2}{1+(k+2)(k+1)\alpha^2}p_{n,k+1},\quad \forall n\geq 0
\quad 1\leq k\leq n+1\\
p_{n,k}=0,\quad \forall k>n,\quad p_{0,0}=1
\end{array}
\right.
\end{align*}

furthermore

\begin{align*}
\left\{
\begin{array}{ll}
p_{n,k}=d_kT^{\left(\frac{n-k}{2}\right)}_k,\quad n\geq 0,\quad 0\leq k\leq n,\quad n+k\quad even\\
p_{n,k}=0,\quad otherwise\\
d_k=\prod_{i=0}^{k}\frac{1}{1+i(i-1)\alpha^2}\\
T^{(h)}_k=\sum_{i=0}^{k}\frac{d_{i+1}-d_{i+2}}{d_i}T^{(h-1)}_{i+1},\quad \forall k\geq 0,h\geq 1;\quad T^{(0)}_k=1,\quad \forall k\geq 0
\end{array}
\right.
\end{align*}

\end{proposition}

\subsection{Steady state}

\noindent
Taking into account that a $M/M/\infty$ queue has an exponential rate of convergence to equilibrium, see \cite{Joulin}, for applications it is interesting to know the behavior of the degree distribution in the steady state.
We know that in the steady state the number of the nodes follows the Poisson distribution with parameter $\frac{\lambda}{\mu}$, then
we can recover the degree distribution of the generic node starting from the events $N=l>0$.\\
\noindent
Let us consider the steady state and label the nodes with the natural sequence choosing the index $n$ for node position. From Eq.(\ref{equ:outdegreebdp}), conditioning to $N=l$, we have that the distribution in the steady state verifies

\begin{equation}\label{equ:statcondbdp}
0=\pi^n_{k-1,l}\frac{\lambda}{k}+ \pi^n_{k+1,l}\mu (k+1)+\pi^n_{k,l}\left(-\frac{\lambda}{1+k}-\mu k\right)
\end{equation}

\noindent
Moreover if $k\geq l-n\geq 0$ we have $\pi^n_{k,l}=0$ because the out-degree of $n$-th node can not be greater than $l-n$
conditioned to the event $N=l$.\\

\begin{proposition}

If $\pi^n_{k,l}$ satisfies Eq.(\ref{equ:statcondbdp}) then for $k<l-n$ and $n<l$ with $l,n=1,2,3,...$ it results

\begin{align*}
\left\{
\begin{array}{ll}
\pi^n_{k,l}=\pi_0\left(\frac{\lambda}{\mu}\right)^{k}\frac{1}{(k!)^2}\\
\pi^n_{0,l}=\left(\sum_{k=1}^{l-n}\left(\frac{\lambda}{\mu}\right)^{k}\frac{1}{(k!)^2}\right)^{-1}
\end{array}
\right.
\end{align*}

furthermore for the unconditional distribution $\pi^n_k$ it results

\begin{align*}
&\pi^n_{k}=\left(\frac{\lambda}{\mu}\right)^{k}\frac{e^{-\frac{\lambda}{\mu}}}{(k!)^2}
\sum_{l=n+k}^{+\infty}\frac{(\frac{\lambda}{\mu})^l}{l!\sum_{m=1}^{l-n}(\frac{\lambda}{\mu})^{m}\frac{1}{(m!)^2}}
\end{align*}

\end{proposition}

\subsection{Recover connections from the root}

\noindent
In the general case it is possible for a node to have in-degree $0$ at a generic time $t$. This generates a temporary inefficiency of the streaming service.
In this case we create with probability $1$ a connection from the root to this node and restore the service. We will not modify
the out-degree process for the root, we consider the act of connection recovery only as an extra work for the root.\\
Let us consider the steady state. We know that the number $N$ of the nodes (\emph{except the root}) is distributed in accordance with Poisson distribution with parameter $\frac{\lambda}{\mu}$. Moreover the number of connections that each node established to benefit of the stream service is
related to transitions of the out-degree processes of the nodes that was early in the network. Because the out-degree processes are independent these numbers of (in) connections is a sum of Bernoullian random variables.\\
Let $X^{-}_i$ the random number of in-arcs of the generic node, from the fact that all nodes have the same distribution for the life time and from the minimum order statistic, we deduce that the extra work for the root is:

\begin{equation}\label{extraw}
W_{root}=\sum_{i=1}^{N}\frac{1}{X^-_i}
\end{equation}

\noindent
By virtue of Anscombe transform and of a generalized Hoeffding inequality \cite{Hoeffding_gen} we can established the following result.

\begin{proposition}

Let $W_{root}$ as in Eq.(\ref{extraw})  then

\begin{align*}
&E(W_{root})\sim \sqrt{\frac{\lambda}{\mu}}
\end{align*}

\end{proposition}

\bigskip
\noindent
\textbf{Acknowledgement}
I am grateful to my Ph.D. advisor Marco Isopi to encourage me for all time.
I wish also to thank Lorenzo Favalli of Electronic Department at University of Pavia, for his support on understanding content delivery networks characteristic and many explanations about peer-to-peer media streaming.

\end{document}